\documentclass[twoside]{article} 
\usepackage{graphicx}
\date{} 
\title{The asymptotic expansion of Kr\"atzel's integral and an integral related to an extension of the Whittaker function}
\author{\sc R. B.\ Paris \\
{\em Division of Computing and Mathematics,} \\
{\em Abertay University, Dundee DD1 1HG, UK}}
\begin{document}
\def\f#1#2{\mbox{${\textstyle \frac{#1}{#2}}$}}
\def\dfrac#1#2{\displaystyle{\frac{#1}{#2}}}
\def\boldal{\mbox{\boldmath $\alpha$}}
\newcommand{\bee}{\begin{equation}}
\newcommand{\ee}{\end{equation}}
\newcommand{\lam}{\lambda}
\newcommand{\ka}{\kappa}
\newcommand{\al}{\alpha}
\newcommand{\ba}{\beta}
\newcommand{\la}{\lambda}
\newcommand{\ga}{\gamma}
\newcommand{\eps}{\epsilon}
\newcommand{\fr}{\frac{1}{2}}
\newcommand{\fs}{\f{1}{2}}
\newcommand{\g}{\Gamma}
\newcommand{\br}{\biggr}
\newcommand{\bl}{\biggl}
\newcommand{\ra}{\rightarrow}
\newcommand{\gtwid}{\raisebox{-.8ex}{\mbox{$\stackrel{\textstyle >}{\sim}$}}}
\newcommand{\ltwid}{\raisebox{-.8ex}{\mbox{$\stackrel{\textstyle <}{\sim}$}}}
\renewcommand{\topfraction}{0.9}
\renewcommand{\bottomfraction}{0.9}
\renewcommand{\textfraction}{0.05}
\newcommand{\mcol}{\multicolumn}
\date{}
\maketitle
\pagestyle{myheadings}
\markboth{\hfill \sc R. B.\ Paris  \hfill}
{\hfill \sc Kr\"atzel's integral\hfill}
\begin{abstract}
We consider the asymptotic expansion of  Kr\"atzel's integral
\[F_{p,\nu}(x)=\int_0^\infty t^{\nu-1} e^{-t^p-x/t}\,dt\qquad (|\arg\,x|<\pi/2),\]
for $p>0$ as $|x|\to \infty$ in the sector $|\arg\,x|<\pi/2$ employing the method of steepest descents. An alternative derivation of this expansion is given using a Mellin-Barnes integral approach. The cases $p<0$, $\Re (\nu)<0$  and when $x$ and $\nu$ ($p>0$) are both large are also considered.

A second section discusses the asymptotic expansion of an integral involving a modified Bessel function that has recently been introduced as an 
extension of the Whittaker function $M_{\ka,\mu}(z)$.
Numerical examples are provided to illustrate the accuracy of the various expansions obtained.

\vspace{0.3cm}

\noindent {\bf Mathematics subject classification (2010):} 30E15, 33E20, 33C15, 34E05, 41A60
\vspace{0.1cm}
 
\noindent {\bf Keywords:} Kr\"atzel's integral, asymptotic expansions, steepest descents, Mellin-Barnes integral, modified Bessel function, extended Whittaker function
\end{abstract}

\vspace{0.3cm}

\noindent $\,$\hrulefill $\,$

\vspace{0.3cm}

\begin{center}
{\bf 1.\ Introduction}
\end{center}
\setcounter{section}{1}
\setcounter{equation}{0}
\renewcommand{\theequation}{\arabic{section}.\arabic{equation}}
The first part of this paper is concerned with the asymptotic expansion of Kr\"atzel's integral defined by\footnote{In \cite{KST} this integral was denoted by $Z_p^\nu(x)$.}
\bee\label{e11}
F_{p,\nu}(x)=\int_0^\infty t^{\nu-1} \exp\,\bl(-t^p-\frac{x}{t}\br)\,dt\qquad (|\arg\,x|<\fs\pi),
\ee
where $p\in {\bf R}$, $\nu\in {\bf C}$, with $\Re (\nu)<0$ when $p\leq 0$. For $p\geq 1$, the function $F_{p,\nu}(x)$ was introduced by Kr\"atzel \cite{K} as a kernel of the integral transform
\[\int_0^\infty F_{p,\nu}(xt) \phi(t)\,dt \qquad (x>0).\]
He derived the inversion formula for this transform and constructed its convolution to give applications to the solution of some ordinary differential equations. Further investigations of $F_{p,\nu}(x)$ were carried out by Kilbas {\em et al.} \cite{KST} who established its connection with the Wright function ${}_1\Psi_0(a,b;z)$ defined in Section 3.

Particular cases of (\ref{e11}) are: when $p=1$, $x=z^2/4$ 
\[F_{1,\nu}(z^2/4)=2(\fs z)^\nu K_\nu(z),\]
where $K_\nu(z)$ is the modified Bessel function (or Macdonald function); and when $p=0$ and $p=-1$
\bee\label{e12}
F_{0,\nu}(x)=x^\nu \g(-\nu),\qquad F_{-1,\nu}(x)=(1+x)^\nu \g(-\nu) \qquad (\Re (\nu)<0).
\ee 
We determine the asymptotic expansion of (\ref{e11}) when $p>0$ for $|x|\to\infty$ in the sector $|\arg\,x|<\fs\pi$ by application of the method of steepest descents, with the coefficients in the expansion determined by Lagrange inversion. We present an alternative derivation of the expansion when $p>0$ using a Mellin-Barnes integral representation of $F_{p,\nu}(x)$  given in \cite{KST}. The case when $p<0$ is discussed in Section 4. We also consider the expansion in the case when the parameter $\nu$ is large of $O(x^{p/(p+1)})$ in Section 5.

In the second part, we consider an integral also involving the modified Bessel function $K_\nu(z)$ that arose in an extension of the Whittaker function $M_{\ka,\mu}(z)$ discussed in \cite{DP}. The integral in question is
\[I(x)=\sqrt{\frac{2p}{\pi}} \int_0^1 t^{a-\frac{1}{2}}(1-t)^{b-\frac{1}{2}} e^{-xt} K_\nu\bl(\frac{p}{t(1-t)}\br)\,dt
\qquad(\nu\geq-\fs,\ p>0).\]
The leading approximation of $I(x)$ for large $x>0$ was obtained by application of the saddle-point method in \cite[\S 3]{DP}. Here we determine its asymptotic expansion valid for $|x|\to\infty$ in the sector $|\arg\,x|<\fs\pi$
employing a different approach.

\vspace{0.6cm}

\begin{center}
{\bf 2. The expansion of $F_{p,\nu}(x)$ when $p>0$}
\end{center}
\setcounter{section}{2}
\setcounter{equation}{0}
\renewcommand{\theequation}{\arabic{section}.\arabic{equation}}
Let us consider first the case $p>0$ in (\ref{e11}). We make the change of variable $t=\tau (x/p)^{1/(p+1)}$ to find
\bee\label{e21}
F_{p,\nu}(x)=(x/p)^{\nu/(p+1)} \int_0^\infty \tau^{\nu-1} e^{-X \psi(\tau)}\,d\tau,
\ee
where
\bee\label{e21a}
\psi(\tau):=\frac{\tau^p}{p}+\frac{1}{\tau},\qquad X:=p^{1/(p+1)} x^{p/(p+1)}.
\ee
The integral in (\ref{e21}) converges at both endpoints provided $\Re (X)>0$; that is, when $|\arg\,x|<(p+1)\pi/(2p)$
and so represents the analytical continuation of (\ref{e11}) into this wider sector.

Saddle points of the exponential factor occur when $\psi'(\tau)=0$; that is, when $\tau^{p-1}-1/\tau^2=0$ whence
\[\tau_{sk}=e^{2\pi ik/(p+1)}\qquad (k=0, \pm 1, \pm 2, \ldots  ).\]
The saddle points are consequently situated on the unit circle in the $\tau$-plane, with the saddles corresponding to higher $k$ values being situated on different Riemann surfaces. The only contributory saddle for $F_{p,\nu}(x)$ is found to be that corresponding to $k=0$; see Fig.~1 for typical examples of the steepest descent path from the origin passing to infinity in the direction $\arg\,\tau=-\theta/(p+1)$, where $\theta=\arg\,x$. 

\begin{figure}[th]
	\begin{center}	{\tiny($a$)}\includegraphics[width=0.35\textwidth]{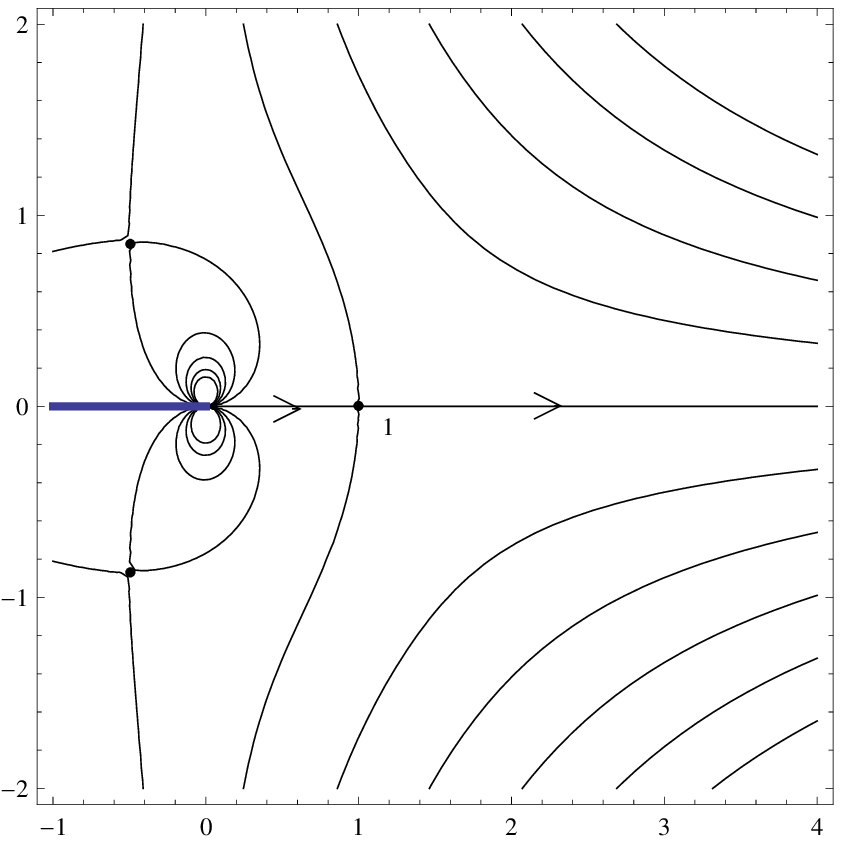}
	\qquad
	{\tiny($b$)}\includegraphics[width=0.35\textwidth]{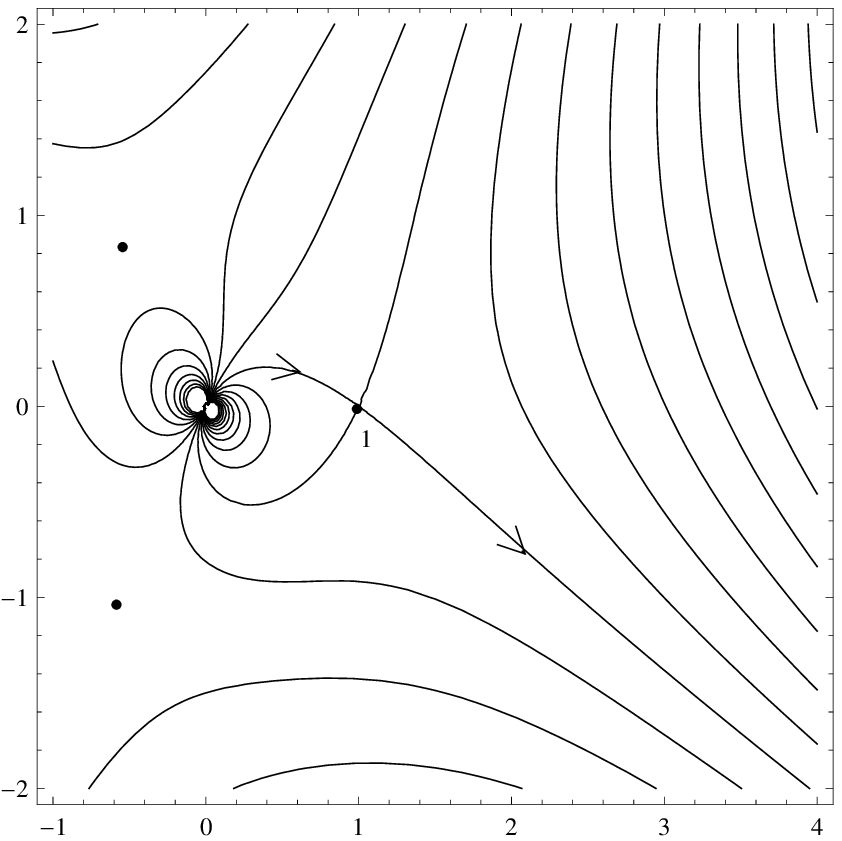} 
\caption{\small{The path of steepest descent from the origin passing over the saddle at $\tau=1$ when $p=2$:
($a$) $\arg\,X=0$ and ($b$) $\arg\,X=\pi/3$. The arrows denote the direction of integration and the dots represent the saddle points situated at $\tau=1, e^{\pm 2\pi i/3}$.}}
\end{center}
\end{figure}

Noting that $\psi(1)=(p+1)/p$, we introduce the new variable $w$ by
\[
\frac{1}{2}w^2=\psi(\tau)-\psi(1)=\frac{\psi''(1)}{2!} (\tau-1)^2+\frac{\psi'''(1)}{3!} (\tau-1)^3+\cdots
\]
so that
\[F_{p,\nu}(x)=(x/p)^{\nu/(p+1)} e^{-X\psi(1)} \int_{-\infty}^\infty e^{-X\tau^2/2} \tau^{\nu-1} \frac{d\tau}{dw}\,dw.\]
Inversion of the above expansion using the {\tt InverseSeries} command in {\em Mathematica} (essentially Lagrange inversion) yields
\bee\label{e22}
\tau=1+\frac{1}{(p+1)^{1/2}} \sum_{k=1}^\infty A_k w^k,
\ee
where the first few coefficients $A_k$ are
\[A_1=1,\quad A_2=\frac{4-p}{6(p+1)^{1/2}},\quad A_2=\frac{26-19p+2p^2}{72(p+1)},\]
\[A_3=\frac{2}{135(p+1)^{3/2}} (46-66p+21p^2-p^3),\]
\[A_4=\frac{1}{3456(p+1)^2} (1252-2956p+1881p^2-316p^3+4p^4).\]
Then we have
\bee\label{e23}
\frac{d\tau}{dw}=\frac{1}{(p+1)^{1/2}} \sum_{k=1}^\infty kA_k w^{k-1}.
\ee

Combining (\ref{e22}) and (\ref{e23}), we obtain
\[\tau^{\nu-1}\frac{d\tau}{dw}\stackrel{e}{=} \frac{1}{(p+1)^{1/2}} \sum_{k=0}^\infty B_k w^{2k},\]
where $\stackrel{e}{=}$ denotes the inclusion of only the even powers of $w$ (since the odd powers will not contribute to the integral). The first few coefficients $B_k$ are given by:
\[B_0=1,\quad B_1=\frac{1}{24(p+1)} (2+12\nu^2-12\nu(p-1)-7p+2p^2),\]
\[B_2=\frac{1}{3456(p+1)^2} \bl\{4+144\nu^4-480\nu^3(p-1)-172p+417p^2-172p^3+4p^4\]
\[\hspace{5cm}+120\nu^2(4-11p+4p^2)-24\nu(-6+41p(1-p)+6p^3)\br\}.\]
It is evident that the rapidly increasing complexity of the higher coefficients prevents their explicit representation. However, when dealing with specific values of the parameters $p$ and $\nu$, it is quite feasible to generate numerically many coefficients $B_k$ using the numerical inversion procedure described above.

Then
\[F_{p,\nu}(x)\sim (x/p)^{\nu/(p+1)} e^{-X\psi(1)}\sum_{k=0}^\infty B_k \int_{-\infty}^\infty e^{-Xw^2/2}w^{2k}dw\]
\bee\label{e24}
=\sqrt{\frac{2\pi}{(p+1)X}}\,\frac{e^{-X\psi(1)}}{(p/x)^{\nu/(p+1)}} \sum_{k=0}^\infty \frac{(\fs)_k B_k}{(X/2)^k}\hspace{0.3cm}
\ee
as $x\to+\infty$, where $X$ is defined in (\ref{e21a}) and $(a)_k=\g(a+k)/\g(a)$ is Pochhammer's symbol. We observe that the argument of the exponential factor is
\[X\psi(1)=(p+1)\,(x/p)^{p/(p+1)}.\]
The contributory saddle continues to be $\tau_{s0}=1$ when $x$ is complex satisfying $|\arg\,x|<\fs\pi$; see Fig.~1(b). Hence the expansion in (\ref{e24}) continues to hold for $|x|\to\infty$ in the sector $|\arg\,x|<\fs\pi$.

To illustrate the accuracy of the expansion (\ref{e24}) we present the cases $p=2$, $\nu=\fs$ and $p=3$, $\nu=\f{3}{2}$.
In Table 1 we show the values of the coefficients $B_k$ for $k\leq 6$ obtained by numerical inversion. In Table 2
we give the absolute relative error in the computation of $F_{p,\nu}(x)$ using (\ref{e24}) for different truncation index $k$. The exact value of $F_{p,\nu}(x)$ was determined by high-precision evaluation of the integral in (\ref{e11}).

\begin{table}[bh]
\caption{\footnotesize{The coefficients $B_k$ (with $B_0=1$) for $1\leq k\leq 6$.}}
\begin{center}
\begin{tabular}{|c|l|l|}
\hline
&&\\[-0.3cm]
\mcol{1}{|c|}{$k$} & \mcol{1}{c|}{$p=0.75,\ \nu=1.50$} & \mcol{1}{c|}{$p=1.50,\ \nu=0.50$} \\
[.1cm]\hline
&&\\[-0.25cm]
1 & $+6.99404762\times 10^{-1}$ & $-6.66666667\times 10^{-2}$\\
2 & $+4.53883811\times 10^{-2}$ & $+9.07407407\times 10^{-3}$\\
3 & $-1.50268380\times 10^{-3}$ & $-1.24847737\times 10^{-3}$\\
4 & $+9.39269762\times 10^{-5}$ & $+1.63694787\times 10^{-4}$\\
5 & $-6.82774001\times 10^{-6}$ & $-1.97786516\times 10^{-5}$\\
6 & $+5.11051949\times 10^{-7}$ & $+2.08358000\times 10^{-6}$\\
[.1cm]\hline
\end{tabular}
\end{center}
\end{table}

\begin{table}[bh]
\caption{\footnotesize{Values of the absolute relative error in $F_{p,\nu}(x)$ using (\ref{e24}) for different truncation index $k$ when $x=30$. The final row gives the value of $F_{p,\nu}(x)$.}}
\begin{center}
\begin{tabular}{|c|c|c|}
\hline
&&\\[-0.3cm]
\mcol{1}{|c|}{$k$} & \mcol{1}{c|}{$p=0.75,\ \nu=1.50$} & \mcol{1}{c|}{$p=1.50,\ \nu=0.50$} \\
[.1cm]\hline
&&\\[-0.25cm]
0 & $1.679\times 10^{-2}$ & $7.106\times 10^{-3}$\\
1 & $8.182\times 10^{-3}$ & $3.115\times 10^{-4}$\\
2 & $3.478\times 10^{-4}$ & $2.313\times 10^{-5}$\\
3 & $3.966\times 10^{-5}$ & $2.305\times 10^{-6}$\\
4 & $6.844\times 10^{-6}$ & $2.744\times 10^{-7}$\\
5 & $1.504\times 10^{-6}$ & $3.542\times 10^{-8}$\\
[.1cm]\hline
&&\\[-0.25cm]
$F_{p,\nu}(x)$ & $5.70691385\times 10^{-3}$ & $2.67529872\times 10^{-7}$\\
[.1cm]\hline
\end{tabular}
\end{center}
\end{table}

\vspace{0.6cm}

\begin{center}
{\bf 3.\ An alternative derivation of the expansion for $F_{p,\nu}(x)$ when $p>0$}
\end{center}
\setcounter{section}{3}
\setcounter{equation}{0}
\renewcommand{\theequation}{\arabic{section}.\arabic{equation}}
In Kilbas {\em et al.} \cite{KST} it was shown that
\bee\label{e31}
F_{p,\nu}(x)=\frac{1}{2\pi ip} \int_{c-\infty i}^{c+\infty i} \g(s) \g\bl(\frac{s+\nu}{p}\br) x^{-s} ds\qquad c>\max\{0,-\Re (\nu)\}
\ee
valid in the sector $|\arg\,x|<(p+1)\pi/(2p)$. The integration path lies to the right of the poles of $\g(s)$ and $\g((s+\nu)/p)$ situated at $s=-k$ and $s=-kp-\nu$ ($k=0, 1, 2, \ldots $), respectively. This result can be easily established by making use of the Cahen-Mellin integral (see, for example, \cite[p.~89]{PK})
\bee\label{e31a}
e^{-z}=\frac{1}{2\pi i}\int_{c-\infty i}^{c+\infty i} \g(s) z^{-s} ds\qquad (c>0;\ |\arg\,z|<\fs\pi)
\ee
with $z=x/t$ followed by substitution in (\ref{e11}) and evaluation of the resulting integral over $t$ as a gamma function.

Displacement of the integration path to the left over the poles produces (provided $kp+\nu$ does not take on non-negative integer values) \cite{KST}
\[F_{p,\nu}(x)=\frac{1}{p}\sum_{k=0}^\infty \g\bl(\frac{\nu-k}{p}\br)\,\frac{(-x)^k}{k!}+x^\nu\sum_{k=0}^\infty \g(-\nu-kp) \,\frac{(-x^p)^k}{k!}\]
\bee\label{e32}
\ \ =\frac{1}{p}{}_1\Psi_0(-1/p,\nu/p;-x)+x^\nu {}_1\Psi_0(-p,-\nu:-x^p),\ee
where the Wright function ${}_1\Psi_0(a,b;z)$ is defined by
\[{}_1\Psi_0(a,b;z):=\sum_{k=0}^\infty \g(ak+b)\,\frac{z^k}{k!}\qquad (0<a<1).\]
When $kp+\nu$ assumes a non-negative integer value, double poles arise with the result that logarithmic terms will be present.

To determine the asymptotic expansion of the integral in (\ref{e31}) we may displace the integration path as far to the right as we please, since there are no poles of the integrand on the right of this path. On the displaced path, which we denote by $L$, it follows that $|s|$ is everywhere large on $L$. If we define the parameters
\[\ka=\frac{p+1}{p},\quad h=p^{1/p}, \quad \vartheta=\frac{\nu}{p}-\frac{1}{2},\quad {\cal A}_0=(2\pi)^{-1/2} \ka^{\frac{1}{2}-\vartheta} p^{-\vartheta},\]
we have the inverse factorial expansion \cite[p.~39, Lemma 2.2]{PK}
\bee\label{e33}
\g(s)\g\bl(\frac{s+\nu}{p}\br)=2\pi{\cal A}_0(h\ka^\ka)^{-s}\bl\{\sum_{j=0}^{M-1} (-)^j c_j \g(\ka s+\vartheta-j)+\rho_M(s) \g(\ka s+\vartheta-M)\br\},
\ee
where $M$ is a positive integer and $\rho_M(s)=O(1)$ as $|s|\to\infty$ in $|\arg\,s|<\pi$. The coefficient $c_0=1$ and an algorithm for the determination of the higher coefficients $c_j$ ($j\geq 1$) is given in \cite[\S 2.2.4]{PK}, \cite[pp.~321--323]{P}.

Substitution of (\ref{e33}) into (\ref{e31}) taken over the displaced integration path $L$ then yields
\[F_{p,\nu}(x)=\frac{2\pi{\cal A}_0}{p}\bl\{\sum_{j=0}^{M-1}(-)^j c_j\,\frac{1}{2\pi i}\int_L (\ka X)^{-\ka s} \g(\ka s+\vartheta-j)\,ds+R_M(x)\br\},\]
where $X=(hx)^{1/\ka}=p^{1/(p+1)} x^{p/(p+1)}$ (see (\ref{e21a})) and the remainder term $R_M(x)$ is
\[R_M(x)=\frac{1}{2\pi i}\int_L \rho_M(s)(\ka X)^{-\ka s} \g(\ka s+\vartheta-M)\,ds.\]
Application of Lemma 2.7 in \cite[p.~71]{PK} shows that 
\[R_M(x)=O(X^{\vartheta-M}e^{-\ka X}) \qquad (|X|\to\infty,\ |\arg\,X|<\fs\pi).\] 
Hence, upon use of the Cahen-Mellin integral in (\ref{e31a}), we obtain the expansion
\bee\label{e35}
F_{p,\nu}(x)=\frac{2\pi {\cal A}_0}{\ka p}\,X^\vartheta e^{-\ka X} \bl\{\sum_{j=0}^{M-1} (-)^j c_j (\ka X)^{-j}+O(X^{-M})\br\}
\ee
as $|X|\to\infty$ in the sector $|\arg\,X|<\fs\pi$. We note that this sector includes the sector $|\arg\,x|<\fs\pi$.

Routine algebra shows that the leading term of this expansion agrees with the leading term of (\ref{e24}). Consequently we must have the coefficients $c_j$ related to the $B_j$ in (\ref{e24}) given by
\[c_j=(-)^j (2\ka)^j (\fs)_j B_j.\]

\vspace{0.6cm}

\begin{center}
{\bf 4.\ The expansion for $F_{p,\nu}(x)$ when $p<0$}
\end{center}
\setcounter{section}{4}
\setcounter{equation}{0}
\renewcommand{\theequation}{\arabic{section}.\arabic{equation}}
When $p<0$ the expansions of $F_{p,\nu}(x)$ are different according as $-1\leq p<0$ and $p\leq -1$; the situations pertaining to $p=0$ and $p=-1$ being given in (\ref{e12}). In both ranges of the parameter $p$ we require $\Re (\nu)<0$ for convergence of the integral in (\ref{e11}) at the upper limit. 

From \cite[(3.13)]{KST}, we have when $p<0$
\[F_{p,\nu}(x)=-\frac{1}{2\pi ip} \int_{c-\infty i}^{c+\infty i} \g(s) \g\bl(\frac{s+\nu}{p}\br)x^{-s}\,ds\qquad 0<c<-\Re (\nu),\]
where the integration path now separates the poles of $\g(s)$ at $s=-k$ from those of $\g((s+\nu)/p)$ at $s=-kp-\nu$, ($k=0, 1, 2, \ldots$). When $-1\leq p<0$, displacement of the path to the right over the poles at $s=-kp-\nu$ produces
\cite[Theorem 4.2]{KST}
\bee\label{e41}
F_{p,\nu}(x)=x^\nu \sum_{k=0}^\infty \frac{\g(-kp-\nu)}{k!}\,(-x^p)^k \qquad (\Re (\nu)<0,\ -1\leq p<0).
\ee
This convergent series also yields the asymptotic expansion for $|x|\to\infty$ valid in $|\arg\,x|<(p+1)\pi/(2p)$. Hence we have
\[F_{p,\nu}(x)=x^\nu\bl\{\g(-\nu)-\frac{\g(-p-\nu)}{x^{-p}}+O(1/x^{-2p})\br\}\qquad (x\to\infty).\]

When $p\leq -1$, displacement of the path to the left over the poles at $s=-k$ produces \cite[Theorem 4.3]{KST}
\bee\label{e42}
F_{p,\nu}(x)=-\frac{1}{p}\sum_{k=0}^\infty\g\bl(\frac{\nu-k}{p}\br)\,\frac{(-x)^k}{k!}\qquad (\Re (\nu)<0,\ p\leq -1).
\ee
The asymptotic expansion in this case is obtained by displacing the integration path to the right over the poles at $s=-kp-\nu$ to yield 
\bee\label{e43}
F_{p,\nu}(x)\sim x^\nu \sum_{k=0}^\infty \g(-kp-\nu)\, \frac{(-x^p)^k}{k!}\qquad (\Re (\nu)<0,\ p\leq -1)
\ee
as $|x|\to\infty$ in $|\arg\,x|<(p+1)\pi/(2p)$.

It is easily verified that when $p=-1$ the expansions (\ref{e41}) and (\ref{e42}) both reduce to
\[F_{-1,\nu}(x)=(1+x)^\nu \g(-\nu)\]
as stated in (\ref{e12}).

\vspace{0.6cm}

\begin{center}
{\bf 5.\ The expansion of $F_{p,\nu}(x)$ when $\nu$ and $x$ are large}
\end{center}
\setcounter{section}{5}
\setcounter{equation}{0}
\renewcommand{\theequation}{\arabic{section}.\arabic{equation}}
Let $\nu=1+aX$, where $a>0$ is fixed and $X=p^{1/(p+1)} x^{p/(p+1)}$ with $p>0$. Then from (\ref{e21}) we find
\bee\label{e51}
F_{p,\nu}(x)=(x/p)^{\nu/(p+1)} \int_0^\infty \tau^{aX} e^{-X(\tau^p/p+\tau^{-1})}d\tau=(x/p)^{\nu/(p+1)} \int_0^\infty e^{-X\phi(\tau)}d\tau,
\ee
where 
\[\phi(\tau)=\frac{\tau^p}{p}+\frac{1}{\tau}-a\log\,\tau.\]
Saddle points of the exponential factor are given by solutions of the equation
\bee\label{e52}
\tau^{p+1}=1+a\tau.
\ee

Examination of the steepest paths shows that the path from the origin through the contributory saddle $\tau_s>1$ on the positive real axis is similar to that shown in Fig.~1.
Making the change of variable
\[\frac{1}{2}w^2=\phi(\tau)-\phi(\tau_s),\]
we find upon inversion as described in Section 2, followed by differentiation,
\[\frac{d\tau}{dw} \stackrel{e}{=} \frac{1}{\sqrt{\phi_s''}}\,\sum_{k\geq0} C_k w^{2k},\]
where we use the abbreviation $\phi_s=\phi(\tau_s)$ and
\[C_0=1,\qquad C_1=\frac{1}{24\phi_s''^3}  (5\phi_s'''^2-3\phi_s'' \phi_s^{iv}),\]
\[C_2=\frac{1}{3456\phi_s''^6}(385\phi_s'''^4-630\phi_s''\phi_s'''^2\phi_s^{iv}+105\phi_s''^2(\phi_s^{iv})^2+168\phi_s''^2\phi_s'''\phi_s^v-24\phi_s''^3\phi_s^{vi}).\]
Then
\[F_{p,\nu}(x)=(x/p)^{\nu/(p+1)}\,\frac{e^{-X\phi(\tau_s)}}{\sqrt{\phi_s''}} \int_{-\infty}^\infty e^{-Xw^2/2} \,\frac{d\tau}{dw}\,dw\]
\[\sim (x/p)^{\nu/(p+1)}\,e^{-X\phi(\tau_s)}\,\sqrt{\frac{2\pi}{X\phi_s''}} \sum_{k\geq0}\frac{(\fs)_k C_k}{(X/2)^k}\]
as $x\to\infty$. 

Upon noting that
\[\phi(\tau_s)=\frac{p+1}{p\tau_s}+\frac{a}{p}-a\log\,\tau_s,\qquad \phi_s''=\frac{1}{\tau_s^2}\bl(\frac{p+1}{\tau_s}+ap\br),\]
where use has been made of (\ref{e52}), we finally obtain the expansion
\bee\label{e53}
F_{p,\nu}(x)\sim \sqrt{\frac{2\pi}{X\phi_s''}}\,\frac{e^{-\ka X/\tau_s-aX/p}}{(p/x)^{\nu/(p+1)}}\,\tau_s^{aX} \sum_{k\geq 0}\frac{(\fs)_k C_k}{(X/2)^k}
\ee
as $x\to\infty$ when $\nu=1+aX$, $X=p^{1/(p+1)} x^{p/(p+1)}$ and $\ka=(p+1)/p$. It is found that the steepest descent path when $x$ is complex is very similar to that shown in Fig.~1(b), with the real positive saddle $\tau_s$ still being the only contributory saddle when $|\arg\,x|<\fs\pi$. Hence the expansion (\ref{e53}) holds for complex $x$ when $|\arg\,x|<\fs\pi$.

Since, in general, the saddle $\tau_s$ has to be obtained by numerical solution of (\ref{e52}) it is reasonable to determine the coefficients $C_k$ numerically using the {\tt InverseSeries} command in {\em Mathematica}. As an example, we display the coefficients in Table 3 for two cases of $(p,a)$ and in Table 4 the corresponding values of the absolute relative error in the expansion (\ref{e53}) for complex $x$.
\begin{table}[bh]
\caption{\footnotesize{The coefficients $C_k$ (with $C_0=1$) for $1\leq k\leq 5$.}}
\begin{center}
\begin{tabular}{|c|l|l|}
\hline
&&\\[-0.3cm]
\mcol{1}{|c|}{$k$} & \mcol{1}{c|}{$p=2,\ a=1$} & \mcol{1}{c|}{$p=3,\ a=0.50$} \\
[.1cm]\hline
&&\\[-0.25cm]
1 & $-5.12469537\times 10^{-2}$ & $-0.09820959\times 10^{-2}$\\
2 & $+4.75849832\times 10^{-3}$ & $+5.38592491\times 10^{-3}$\\
3 & $-6.64820214\times 10^{-5}$ & $+1.51551186\times 10^{-3}$\\
4 & $-4.28119981\times 10^{-5}$ & $-4.26687479\times 10^{-4}$\\
5 & $+5.50086099\times 10^{-6}$ & $+3.10584484\times 10^{-5}$\\
[.1cm]\hline
\end{tabular}
\end{center}
\end{table}

\begin{table}[h]
\caption{\footnotesize{Values of the absolute relative error in $F_{p,\nu}(x)$ when $\nu=1+aX$ using (\ref{e53})  when $x=30e^{i\theta}$ for different $\theta$ and truncation index $k=3$.}}
\begin{center}
\begin{tabular}{|c|c|c|}
\hline
&&\\[-0.3cm]
\mcol{1}{|c|}{} & \mcol{1}{c|}{$p=2,\ a=1$} & \mcol{1}{c|}{$p=3,\ a=0.50$} \\
\mcol{1}{|c|}{$\theta/\pi$} & \mcol{1}{c|}{$\tau_s=1.32471796$} &\mcol{1}{c|}{$\tau_s=1.11734904$} \\
[.1cm]\hline
&&\\[-0.25cm]
0 &    $1.870\times 10^{-7}$ & $5.319\times 10^{-7}$\\
0.10 & $1.874\times 10^{-7}$ & $5.326\times 10^{-7}$\\
0.20 & $1.886\times 10^{-7}$ & $5.349\times 10^{-7}$\\
0.30 & $1.905\times 10^{-7}$ & $5.385\times 10^{-7}$\\
0.40 & $1.931\times 10^{-7}$ & $5.431\times 10^{-7}$\\
0.45 & $1.947\times 10^{-7}$ & $5.457\times 10^{-7}$\\
[.1cm]\hline
\end{tabular}
\end{center}
\end{table}

\vspace{0.6cm}
  
\begin{center}
{\bf 6.\ An integral connected with an extended Whittaker function}
\end{center}
\setcounter{section}{6}
\setcounter{equation}{0}
\renewcommand{\theequation}{\arabic{section}.\arabic{equation}}
In a recent paper \cite{DP}, an extension of the Whittaker function was introduced by
\[M_{\ka,\mu}^{(p,\nu)}(z)=\frac{x^{\mu+\frac{1}{2}} e^{x/2}}{B(\mu-\ka+\frac{1}{2},\mu+\ka+\frac{1}{2})}\sqrt{\frac{2p}{\pi}} \int_0^1 t^{\mu+\ka-1} (1-t)^{\mu-\ka-1} e^{zt} K_{\nu+\frac{1}{2}}\bl(\frac{p}{t(1-t)}\br)\,dt,\]
where $K_\nu(z)$ is the modified Bessel function and $B(a,b)=\g(a)\g(b)/\g(a+b)$ is the beta function. This representation is valid for $p\geq 0$, $\nu\geq 0$ and $2\mu\neq -1, -2, \ldots\ $. In the case $p=\nu=0$, the extended function reduces to the usual Whittaker function $M_{\ka,\mu}(x)$ given by \cite[p.~337]{DLMF}
\[M_{\ka,\mu}(z)=\frac{x^{\mu+\frac{1}{2}} e^{-z/2}}{B(\mu-\ka+\frac{1}{2},\mu+\ka+\frac{1}{2})} \int_0^1 t^{\mu+\ka-\frac{1}{2}} (1-t)^{\mu-\ka-\frac{1}{2}} e^{zt}dt.\]
Our aim is to determine the asymptotic expansion of the related integral
\bee\label{e61}
I(a,b;z)=\sqrt{\frac{2p}{\pi}} \int_0^1 t^{a-\fr} (1-t)^{b-\fr} e^{zt} K_\nu\bl(\frac{p}{t(1-t)}\br)\,dt\qquad (\nu\geq -\fs,\ p>0),
\ee
valid for arbitrary parameters $a$ and $b$,
as $|z|\to \infty$. This integral corresponds to that appearing in the definition of $M_{\ka,\mu}^{(p,\nu)}(x)$ when $a=\mu+\ka-\fs$, $b=\mu-\ka-\fs$ and $\nu\to\nu+\fs$.
By making the change of variable $t\to 1-t$, it is easily seen that $I(a,b;z)$ satisfies the following Kummer-type transformation
\bee\label{e61c}
I(a,b;z)=e^z I(b,a;-z).
\ee

We first consider the expansion of $I(a,b;z)$ for large $z$ in the sector $|\arg (-z)|<\fs\pi$. We set $z=-x$ and note that
for large $x$ the dominant behaviour of the integrand results from a neighbourhood of $t=0$.  From \cite[(10.40.2)]{DLMF} we have the expansion
\[K_\nu\bl(\frac{p}{t(1-t)}\br)\sim \sqrt{\frac{\pi}{2p}} \exp\,\bl[-\frac{p}{t}-\frac{p}{1-t}\br] \sum_{k=0}^\infty \frac{a_k(\nu)}{p^k}\,(t(1-t))^{k+\fr}\]
as $t\to 0$, where $a_0(\nu)=1$ and
\[a_k(\nu)=\prod_{r=1}^k \frac{4\nu^2-(2r-1)^2}{k!\ 2^{3k}}=(-1)^k\frac{(\fs+\nu)_k(\fs-\nu)_k}{k! 2^k}\quad (k\geq 1).\]
Substitution of this expansion into the integral (\ref{e61}) then leads to
\[I(a,b;-x)\sim \sum_{k=0}^\infty \frac{a_k(\nu)}{p^k} \int_0^1 t^{a+k}(1-t)^{b+k} e^{-xt-p/t}\,e^{-p/(1-t)}dt\]
\[\hspace{1.1cm}=\sum_{k=0}^\infty \frac{a_k(\nu)}{p^k x^{a+k+1}} \int_0^x \tau^{a+k} e^{-(\tau+px/\tau)} g_k(\tau/x)\,d\tau,\]
where
\bee\label{e61b}
g_k(w):=(1-w)^{b+k} e^{-p/(1-w)}=e^{-p} (1-w)^{b+k} e^{-pw/(1-w)}.
\ee

We can expand $g_k(w)$ about $w=0$ to obtain
\[g_k(w)=e^{-p} \sum_{r=0}^\infty \frac{(-1)^r c_r(k)}{r!} w^r\qquad (|w|<1),\]
where the first few coefficients $c_r(k)$ are given by
\[c_0(k)=1,\qquad c_1(k)=p+\ba, \]
\[c_2(k)=(p+\ba)^2-(p+\ba)-p,\]
\bee\label{e61a}
c_3(k)=(p+\ba)^3-6(p+\ba)^2+(3b+2)(p+\ba)+4p,
\ee
with $\ba:=b+k$. Then
\[I(a,b;-x)\sim e^{-p}\sum_{k=0}^\infty \frac{a_k(\nu)}{p^k x^{a+k+1}} \sum_{r=0}^\infty \frac{(-1)^r c_r(k)}{r! \,x^r} \int_0^x \tau^{a+k+r} e^{-(\tau+px/\tau)}d\tau.\]

We now extend the upper limit $x$ to $\infty$, thereby introducing an error of $O(e^{-x})$ when $|\arg\,x|<\fs\pi$. Employing the integral representation of the modified Bessel function \cite[(10.32.10)]{DLMF}
\[2(\fs z)^{-\nu}K_\nu(z)=\int_0^\infty e^{-(\tau+z^2/(4\tau))}\,\frac{d\tau}{\tau^{\nu+1}}\qquad (|\arg\,z|<\f{1}{4}\pi)\]
and the fact that $K_{-\nu}(z)=K_\nu(z)$, we find (with $\omega:=a+1+k+r$)
\[I(a,b;-x)\sim 2e^{-p} \sum_{k=0}^\infty \frac{a_k(\nu)}{p^k x^{a+k+1}} \sum_{r=0}^\infty \frac{(-1)^r c_r(k)}{r!\,x^r}\,(px)^{\mu/2} K_\omega(2\sqrt{px})\]
\bee\label{e62}
\hspace{1.2cm}=2e^{-p} \bl(\frac{p}{x}\br)^{(a+1)/2} \sum_{k=0}^\infty a_k(\nu)\sum_{r=0}^\infty \frac{(-1)^r c_r(k) p^r}{r!\,(px)^{(k+r)/2}}\,K_\omega(2\sqrt{px}).
\ee
\begin{table}[th]
\caption{\footnotesize{The absolute relative error in the computation of $I(a,b;-x)$ from (\ref{e61}) for different truncation index $r$ in the asymptotic expansion (\ref{e63}) when $a=b=p=1$, $\nu=4/3$.}}
\begin{center}
\begin{tabular}{|c|c|c|c|}
\hline
&&&\\[-0.3cm]
\mcol{1}{|c|}{$r$} & \mcol{1}{c|}{$x=100$} & \mcol{1}{c|}{$x=200$}& \mcol{1}{c|}{$x=500$} \\
[.1cm]\hline
&&&\\[-0.25cm]
0 & $1.912\times 10^{-1}$ & $1.181\times 10^{-1}$ & $6.658\times 10^{-2}$ \\
1 & $2.489\times 10^{-2}$ & $1.147\times 10^{-2}$ & $4.267\times 10^{-3}$ \\
2 & $5.124\times 10^{-3}$ & $1.484\times 10^{-3}$ & $3.167\times 10^{-4}$ \\
3 & $1.732\times 10^{-4}$ & $3.503\times 10^{-5}$ & $4.647\times 10^{-6}$ \\
4 & $4.387\times 10^{-5}$ & $5.780\times 10^{-6}$ & $4.499\times 10^{-7}$ \\
5 & $1.236\times 10^{-5}$ & $1.065\times 10^{-6}$ & $4.876\times 10^{-8}$ \\
6 & $4.694\times 10^{-6}$ & $2.620\times 10^{-7}$ & $6.977\times 10^{-9}$ \\
[.1cm]\hline
\end{tabular}
\end{center}
\end{table}

The double series can be summed ``diagonally'' (cf. \cite[p.~58]{LJS}) to produce the final result
\bee\label{e63}
I(a,b;-x)\sim 2e^{-p}\bl(\frac{p}{x}\br)^{(a+1)/2} \sum_{r=0}^\infty \frac{D_r K_{a+1+r}(2\sqrt{px})}{(px)^{r/2}}
\ee
as $|x|\to\infty$ in $|\arg\,x|<\fs\pi$, where $D_0=1$ and
\[D_r=\sum_{\ell=0}^r \frac{(-1)^\ell p^\ell}{\ell!}\, a_{r-\ell}(\nu) c_\ell(r-\ell)\qquad (r\geq 1).\]
The expansion valid in $|\arg\,z|<\fs\pi$ immediately follows from (\ref{e63}) and the Kummer-type transformation in (\ref{e61c}).

The leading asymptotic form is consequently
\bee
I(a,b;-x)\sim 2e^{-p} \bl(\frac{p}{x}\br)^{(a+1)/2} K_{a+1}(2\sqrt{px})\sim \frac{\sqrt{\pi}p^{a/2+1/4}}{x^{a/2+3/4}}\,e^{-p-2\sqrt{px}}
\ee
as $|x|\to\infty$ in $|\arg\,x|<\fs\pi$.
A result equivalent to the above was obtained by the saddle-point method in \cite[\S 3]{DP}.

In Table 5 we show  an example of the values of the absolute relative error in the computation of $I(a,b;-x)$ from (\ref{e61}) and (\ref{e63}) for different $x$ and truncation index $r$ in the asymptotic sum. The calculation of the coefficients $c_r(k)$ in (\ref{e61b}) can be carried out numerically when specific values of the parameters are known using the {\tt InverseSeries} command in {\it Mathematica}. 

To conclude, we note that an extended confluent hypergeometric function  was introduced in \cite{CQSP} that involved the integral
\[J(a,b;z)=\int_0^1 t^a(1-t)^b e^{zt} \exp \bl(\frac{p}{t(1-t)}\br)\,dt\qquad (p>0).\]
The same arguments employed above produce the expansion
\bee\label{e65}
J(a,b;-x)\sim 2e^{-p}\bl(\frac{p}{x}\br)^{(a+1)/2} \sum_{r=0}^\infty \frac{(-1)^r p^rc_r(0)}{r! (px)^{r/2}}\,K_{a+1+r}(2\sqrt{px})
\ee
as $|x|\to\infty$ in $|\arg\,x|<\fs\pi$, where the coefficients $c_r(0)$ are obtained from (\ref{e61a}) with $\ba=b$.
The expansion valid in $|\arg\,z|<\fs\pi$ follows from (\ref{e65}) and the Kummer-type transformation $J(a,b;z)=e^{z} J(b,a;-z)$.

\end{document}